\documentclass[a4paper,12pt]{amsart}

\usepackage{amsmath,amsfonts,amsthm,amssymb,paralist}
\usepackage[all]{xy}

\newcommand{\Cn}{\mathbb{C}^n}

\newcommand{\C}{\mathbb{C}}
\newcommand{\Z}{\mathbb{Z}}
\newcommand{\CN}{\mathbb{C}^N}

\newcommand{\R}{\mathbb{R}}

\newcommand{\N}{\mathbb{N}}
\newcommand{\Nn}{\mathbb{N}^n}

\newcommand{\dopt}[2]{\frac{\partial #1}{\partial #2}}
\newcommand{\vardop}[3]{\frac{\partial^{|#3|} #1}{\partial {#2}^{#3}}}

\newlength{\extendaxesby}
\DeclareMathOperator{\trord}{tr ord}

\DeclareMathOperator{\jac}{Jac}
\DeclareMathOperator{\ord}{ord}

\DeclareMathOperator{\imag}{Im}

\newtheorem{thm}{Theorem}
\newtheorem{lem}[thm]{Lemma}
\newtheorem{prop}[thm]{Proposition}
\newtheorem{cor}[thm]{Corollary}

\theoremstyle{definition}

\newtheorem{exa}[thm]{Example}

\newtheorem{conj}[thm]{Conjecture}

\def\dbl{[\hskip -1pt [}
\def\dbr{]\hskip -1pt]}
\def \Rk{\text{\rm Rk}\,}

\numberwithin{equation}{section}
 \numberwithin{thm}{section}

\begin{document}

\title{Remarks on the rank properties of formal CR maps}
\author{Bernhard Lamel}
\address{Universit\"at Wien, Fakult\"at f\"ur Mathematik, Nordbergstrasse 15, A-1090 Wien, \"Osterreich}
\email{lamelb@member.ams.org}%
\author{Nordine Mir}
\address{Universit\'e de Rouen, Laboratoire de Math\'ematiques Rapha\"el Salem, UMR 6085 CNRS, Avenue de
l'Universit\'e, B.P. 12, 76801 Saint Etienne du Rouvray, France}
\email{Nordine.Mir@univ-rouen.fr}
\thanks{The first author was supported by the FWF, Projekt P17111}

\begin{abstract}
We prove several new transversality results for formal CR maps between
formal real hypersurfaces in complex space. Both cases of finite and 
infinite type hypersurfaces are tackled in this note.
\end{abstract}

\maketitle
\section*{Introduction}

Given two smooth real hypersurfaces $M,M'$ in $\C^N$ and a
smooth CR mapping $f\colon M\to M'$, it is often of fundamental 
importance to know whether $f$ is
CR-transversal at a given point $p\in M$, that is,
whether the the normal derivative of
the normal component of $f$ at $p$ is non-vanishing. This question naturally
arises from regularity or unique continuation problems for local CR mappings
between smooth or real-analytic hypersurfaces in complex space (see e.g.\ the
works \cite{Forn1, BR6, DF1, H2, BR5, BHR2, HuPa, CR1}).

In this note, we prove several new transversality results for CR mappings
between real hypersurfaces in complex space. Our considerations being purely
formal, we formulate our results in the context of formal CR maps between
formal real hypersurfaces. Two opposite situations are tackled in this note.
Firstly, when the source hypersurface $M$ is of finite type (in the sense of
Kohn \cite{Ko1} and Bloom-Graham \cite{BG1}), we give a new sufficient
condition on $M$ which guarantees that any (formal) holomorphic map sending $M$
into another (formal) real hypersurface in $\C^N$, with non-vanishing Jacobian
determinant, is necessarily CR-transversal (Corollary~\ref{c:jac}). This result
generalizes, in particular, a result due to Baouendi-Rothschild \cite{BR4}.
Secondly, we also deal with the case of infinite type hypersurfaces,
which, to our knowledge, does not seem to have been much studied before (except implicitly in the works of Baouendi-Rothschild \cite{BR7} and Meylan \cite{Meylan95}). 
We prove, among other results, a sharp transversality result for
self-maps of infinite type hypersurfaces (Corollary~\ref{c:transtype2}). We
also provide, in such a setting, a new sufficient condition for a
CR-transversal self-map to be an automorphism
(Theorem~\ref{p:inftypetvimpaut}).

In the next section, we collect the necessary notation and background needed 
in the paper. The results (and proofs) of the results in the finite
type and in the infinite type case can be read independently and are given in
\S \ref{s:finitetype} and \S \ref{s:infinitetype} respectively.

\section{Notation and preliminaries}\label{s:beginning}

For $x=(x_1,\ldots,x_k)\in \C^k$, we denote by $\C \dbl x \dbr$
the ring of formal power series in $x$ and by $\C\{x\}$ the
subring of convergent ones. If $I\subset \C \dbl x\dbr$ is an
ideal and $F:(\C_x^k,0)\to (\C^{k'}_{x'},0)$ is a formal map, 
we define the {\em pushforward} $F_*(I)$ of $I$ to be the ideal in
$\C \dbl x'\dbr$, $x'\in \C^{k'}$, $F_*(I):=\{h\in \C \dbl
x'\dbr:h\circ F\in I\}$. We also define the {\em generic rank} of
$F$, and denote it by $\Rk F$, to be the rank of the Jacobian matrix
${\partial F}/{\partial x}$ regarded as a $\C \dbl x\dbr$-linear
mapping $(\C \dbl x\dbr)^k\to (\C \dbl x\dbr)^{k'}$. Hence $\Rk F$
is the largest integer $r$ such that there is an $r\times r$ minor
of the matrix ${\partial F}/{\partial x}$ which is not  0 as a
formal power series in $x$. Note that if $F$ is convergent, then
$\Rk F$ is the usual generic rank of the map $F$. In addition, for
any complex-valued formal power series $h(x)$, we  denote by $\bar
h (x)$ the formal power series obtained from $h$ by taking complex
conjugates of the coefficients. We also denote by ${\rm ord}\,
h\in \N\cup \{+\infty\}$ the order of $h$ i.e. the smallest
integer $r$ such that $\partial^{\alpha}h(0)=0$ for all $\alpha
\in \N^k$ with $|\alpha|\leq r-1$ and for which
$\partial^{\beta}h(0)\not =0$ for some $\beta \in \N^k$ with
$|\beta|=r$ (if $h\equiv 0$, we simply set ${\rm ord}\, h=+\infty$).
Moreover, if $S=S(x,x')\in \C \dbl x,x'\dbr$, we write ${\rm
ord}_{x}\,S$ to denote the order of $S$ viewed as a power series
in $x$ with coefficients in the ring $\C \dbl x'\dbr$.

For $(Z,\zeta)\in \C^N \times \C^N$, we define the involution
$\sigma : \C\dbl Z,\zeta\dbr \to \C\dbl Z,\zeta\dbr$ by $\sigma
(f)(Z,\zeta):=\bar{f}(\zeta,Z)$.  Let $r\in
\C \dbl Z,\zeta \dbr$ that is {\em 
invariant under the involution} $\sigma$. Such an $r$ is said to
define a {\em formal real hypersurface} through the origin, which
we denote by $M$, if $r(0)=0$ and $r$ has non-vanishing gradient with respect to $Z$ at $0$.
Throughout the note,  
we shall freely write $M\subset \CN$. The complex
space of vectors of $T_0\CN$ which are in the kernel of the
complex linear map $\partial_Zr (0)$ will be denoted by
$T_0^{1,0}M$ (see e.g.\ \cite{BERbook} for the motivation of these definitions).

Throughout this note, it will be convenient to use (formal) {\em
normal coordinates} associated to any formal real hypersurface $M$
of $\CN$ (see e.g.\ \cite{BERbook}). They are
given as follows. There exists a formal change of coordinates in
$\CN\times \CN$ of the form $(Z,\zeta)=(Z(z,w),\bar
Z(\chi,\tau))$, where $Z=Z(z,w)$ is a formal change of coordinates
in $\CN$ and where
$(z,\chi)=(z_1,\ldots,z_{N-1},\chi_1,\ldots,\chi_{N-1})\in \C^{N-1}\times
\C^{N-1}$, $(w,\tau)=(w,\tau)\in \C
\times \C$ so that $M$ is defined through the following defining
equations
\begin{equation}\label{e:coord}
r((z,w),(\chi,\tau))=w-Q(z,\chi,\tau),
\end{equation}
where $Q\in \C \dbl z,\chi,\tau \dbr$ satisfies
\begin{equation}\label{e:normality}
Q(0,\chi,\tau)=Q(z,0,\tau)=\tau.
\end{equation}

In such a coordinate system, we now recall several well-known (invariant) nondegeneracy conditions for the formal real hypersurface $M$. Firstly, $M$ is said to be of {\em finite type} (in the sense of Kohn \cite{Ko1} and Bloom-Graham \cite{BG1}) if $Q(z,\chi,0)\not \equiv 0$. If $M$ is not of finite type, we say that $M$ is of {\em  infinite type}. Following \cite{LM2, LM3}, we say
that $M$ {\em belongs to the class} ${\mathcal C}$ (resp.\ is {\em essentially finite}) if for $k$ large enough the formal $($holomorphic$)$ map $\chi \mapsto
\left(Q_{z^\alpha}(0,\chi,0)\right)_{|\alpha|\leq k}$ is of generic rank $N-1$ (resp.\ is a finite map). Recall also that $M$ is called {\em holomorphically nondegenerate} (in the sense of Stanton \cite{S2}) if the generic rank of the formal map
$(\chi,\tau) \mapsto
\left(Q_{z^\alpha}(0,\chi,\tau)\right)_{|\alpha|\leq k}$ is $N$. 
All these nondegeneracy conditions are known
to be independent of the choice of normal coordinates. Note that if $M\in {\mathcal C}$, then $M$ is necessarily holomorphically nondegenerate. It is also good to point out for the purposes of this paper that formal real hypersurfaces that do not contain any formal curve are examples of hypersurfaces that belong to the class ${\mathcal C}$ (see \cite{LM2, LM3, BERbook} for more details).

If $r,r'\in \left(\C\dbl Z,\zeta\dbr \right)^2$ define two formal real hypersurfaces $M$ and
$M'$, we let ${\mathcal I}(M)$
(resp.\ ${\mathcal I}(M')$) be the ideal generated by $r$ (resp.\
by $r'$). Let $H\colon (\CN,0)\to (\CN,0)$ be a formal
holomorphic map. We associate to the map $H$ another formal map
${\mathcal H}\colon (\CN \times \CN,0)\to (\CN\times \CN,0)$
defined by ${\mathcal H}(Z,\zeta)=(H(Z),\bar{H}(\zeta))$. We say
that {\em $H$ sends $M$ into $M'$} if ${\mathcal I}(M')\subset
{\mathcal H}_*({\mathcal I}(M))$ and write $H(M)\subset M'$. The
formal map $H$ 
is 
{\em CR-transversal} if
\begin{equation}\label{e:crtransverse}
T_0^{1,0}M'+dH(T_0\CN)=T_0\CN,
\end{equation}
where $dH$ denotes the differential of $H$ (at $0$). 

If we choose normal coordinates $Z=(z,w)\in \C^N$ and $Z'=(z',w')\in \C^N$ for
$M$ and $M'$ respectively, and if we write the map $H=(F,G)\in \C^{N-1}\times
\C^N$, then the $H$ is CR-transversal if and only
if $\partial G/\partial w(0)\not \equiv 0$ (see e.g.\ \cite{BR4, ER1}). We refer 
to the last component $G$ of $H$ in any
system of normal coordinates as a {\em transversal $($or normal$)$}
component of $H$.  When $G\equiv 0$, we say that $H$ is {\em transversally
flat}. Following \cite{BR4}, we also say that $H$ is {\em not totally
degenerate} if $z\mapsto \displaystyle \frac{\partial F}{\partial z} (z,0)$ is of generic rank $N-1$.
It is easy to see (and well-known) that the above conditions are independent of
the choice of normal coordinates for $M$, $M'$. Throughout the note, we will
denote by ${\rm Jac}\, H$ the usual Jacobian determinant of $H$.

\section{The finite type case}\label{s:finitetype}

\subsection{A prolongation lemma}
 Our main transversality result in the finite case (Theorem~\ref{t:hypinc} below) is based on the following result. In what follows, for every integer $k$, we denote by $J^k_0(\C^n,\C^d)$ the space of $k$-jets at $0$ of holomorphic maps from $\C^n$ into $\C^d$.
 
\begin{lem}
  \label{l:prolongation}
  Let $A \colon \left( \C^n_z\times\C^m_\chi,(0,0) \right)
  \to \C$ be a formal $($holomorphic$)$ power series, and assume that
  $\ord_z A(z,\chi) = k< \infty$. 
  Then, for every positive integer $d$ and for each 
  $\alpha\in\Nn$, there exists a universal linear map
  \begin{equation}\label{e:prolsol}
    T_{\alpha} \colon J_0^{k+|\alpha|}(\Cn,\C^d) \to 
    \C^d 
  \end{equation}
  $($over the quotient
  field of meromorphic powers series in $\chi)$ such that 
   for every formal power series mapping $b\colon (\C^n\times \C^m, 0)\to \C^d$ and every $v_{\beta}\colon (\C^m,0)\to \C^d$ for $\beta \in \Nn$ satisfying
  \begin{equation}\label{e:prol1}
    \left[\vardop{}{z}{\beta}\left( A(z,\chi)\, b(z,\chi)\right)\right]
    \biggr|_{z=0} = 
    v_\beta (\chi),\quad {\rm for}\, {\rm all}\, \beta\in\Nn,
  \end{equation}
  then
  \begin{equation}
    b_{z^\alpha} (0,\chi) = T_{\alpha} \left( 
    (v_\beta (\chi))_{|\beta|\leq |\alpha| +k} \right).
    \label{e:prolsol2}
  \end{equation}
\end{lem}
\begin{proof}
  Note that the case $\ord_z A(z,\chi) = 0$ is trivial (we can
  just divide by $A$ then), so we assume that
  $k>0$ for the proof.
  
  We order the multiindeces in $\Nn$ 
  by lexicographic ordering; that is, $\alpha < \beta$ if
  and only if $\alpha \not =\beta$ and for the least $j$ with $\alpha_j - \beta_j 
  \neq 0$ this number is negative.

  We first choose $\alpha^0$ to be the minimal
  multiindex $\alpha\in \N^n$ with
  $A_{z^\alpha} (0,\chi)\neq 0$, fix the integer $d$ and prove the lemma by induction on $\ell=|\alpha|$.
  For $\ell=0$, the result is obvious by using \eqref{e:prol1} for $\beta=\alpha^0$.

 Now assume
 that we have constructed the maps $T_{\alpha}$ for
 $|\alpha| < \ell$ for some positive integer $\ell$. 
 It is clear from the chain rule that the equations \eqref{e:prol1} for
 $|\beta| = \ell + k$ only involve the 
 $b_{z^{\gamma}}(0,\chi)$ for $|\gamma| \leq \ell$. We first replace, inside such equations, 
 the  terms $b_{z^\gamma}(0,\chi) $ by $T_\gamma(v_\beta (\chi):|\beta|\leq k+\ell-1)$ for all $\gamma$ of 
 length strictly less than $\ell$. We now 
 compute the coefficient of $b_{z^{\gamma}}(0,\chi)$ for all $\gamma$ of length $\ell$. For every such multiindex, 
 we consider the equation given by \eqref{e:prol1} for $\beta=\alpha^0+\gamma$. Ordering the $\gamma$'s in ascending order, so 
 that $\gamma^1 < \gamma^2 < \dots < \gamma^N$, the above equations yield that for every $j=1,\ldots,N$,
 \[\sum_{\gamma \geq \gamma^j,\atop |\gamma|=\ell}A_{z^{\alpha^0+\gamma-\gamma^j}}(0,\chi)b_{z^\gamma}(0,\chi)=\rho_j\left(v_\beta(\chi):|\beta|\leq \ell+k\right)
 \]
 for some universal linear map
  $\rho_{j} \colon J_0^{k+\ell}(\Cn,\C^d) \to 
    \C^d$  (over ${\rm Frac}\,  \C  \dbl \chi\dbr$)
 The obtained system is 
  triangular in the unknowns $b_{z^{\gamma^j}}(0,\chi)$, $j=1,\ldots,N$ and can be solved. The proof of Lemma~\ref{l:prolongation} is finished.
\end{proof}

\subsection{Statement of the results}

Our main result for finite type hypersurfaces is given by the following.

\begin{thm}\label{t:hypinc} Let $M,M'\subset \C^N$ be formal real hypersurfaces with $M$ belonging to the class ${\mathcal C}$. Then every formal holomorphic map $H\colon(\C^N,0)\to (\C^N,0)$ sending $M$ into $M'$ which is transversally nonflat is necessarily not totally 
  degenerate, and hence CR-transversal.
\end{thm}

As a first consequence of Theorem~\ref{t:hypinc}, we have the following result which provides a partial answer to a question posed in the paper \cite{ER1} (see also Conjecture~\ref{c:conj1} below). 

\begin{cor}\label{c:jac}
 Let $M,M'\subset \C^N$ be as in Theorem~{\rm \ref{t:hypinc}}. Then every formal holomorphic map $H\colon(\C^N,0)\to (\C^N,0)$ sending $M$ into $M'$ satisfying $\jac H\not \equiv 0$ is CR-transversal.
\end{cor}

We also obtain the following transversality result for arbitrary mappings.

\begin{cor}\label{c:newcor}
Let $M,M'\subset \C^N$ be formal real hypersurfaces. Assume that $M$ belongs to the class ${\mathcal C}$ and that 
$M'$ does not contain any formal curve. Then any formal holomorphic map $H\colon (\C^N,0)\to (\C^N,0)$ sending $M$
into $M'$ is either constant or CR-transversal. 
\end{cor}

Theorem~\ref{t:hypinc} improves a transversality result due to Baouendi-Rothschild \cite{BR4} established for essentially finite hypersurfaces. In fact, in this note, we will only show that the map $H$ in Theorem~\ref{t:hypinc}
is not totally degenerate, the CR-transversality of the map being then a consequence of the recent result of Ebenfelt-Rothschild  \cite[Theorem 4.1]{ER1}. We should also point out that when $M$ is assumed to be essentially finite, it was shown in \cite{BR4} that the map is in fact CR-transversal and finite. This need not be the case when $M$ is merely in the class ${\mathcal C}$ as shown by the following example.

\begin{exa}\label{ex:example1}
Let $\psi=(\psi_1,\ldots,\psi_n)\in \left(\C \dbl z\dbr\right)^n$ be a formal map of generic rank $n$, not finite and vanishing at the origin. Consider the formal hypersurface
$M_\psi:=\{(z,w)\in \C^{n}\times \C: {\rm Im}\, w=\sum_{j=1}^{n}|\psi_j(z)|^2\}$ as well as the Heisenberg hypersurface
$\mathbb{H}_n:=\{(z',w')\in \C^n\times \C: {\rm Im}\, w'=\sum_{j=1}^n|z_j'|^2\}$. Then the hypersurface $M_\psi\in {\mathcal C}$ and the formal map $H$ given by
$H(z,w)=(\psi (z),w)$ sends $M_\psi$ into $\mathbb{H}_n$, is CR-transversal, not totally degenerate and obviously not finite.
\end{exa}

In fact, we have the following more precise statement.

\begin{prop}\label{c:add}  Let $M,M'\subset \C^N$ be as in Theorem~{\rm \ref{t:hypinc}} and $H\colon (\C^N,0)\to (\C^N,0)$ be a formal holomorphic map sending $M$ into $M'$.  Then the following are equivalent:
\begin{enumerate}
\item[{\rm (i)}] $H$ is CR transversal;
\item[{\rm (ii)}] $H$ is not totally degenerate.
\end{enumerate}
\end{prop}

Note that Proposition~\ref{c:add} is not true in the more general situation
when $M$ is holomorphically nondegenerate and of finite type (see for instance
the example given in \cite[p.26-27]{ER1}). However, we believe that the
following should be true:

\begin{conj}\label{c:conj1}
 Let $M,M'\subset \C^N$ be formal real hypersurfaces with $M$  holomorphically nondegenerate and of finite type. Then the following are equivalent:
\begin{enumerate}
\item[(i)] $H$ is CR transversal;
\item[(ii)] ${\rm Jac}\, H\not \equiv 0$. 
\end{enumerate}
\end{conj}

We support such a conjecture by proving the following result, which in particular shows that the implication
${\rm (i)}\Rightarrow {\rm (ii)}$ in Conjecture~\ref{c:conj1} is indeed true.

\begin{thm}\label{t:hndhyp}
Let $M,M'\subset \C^N$ be formal real hypersurfaces with $M$ holomorphically nondegenerate.
 Then every formal holomorphic map $H\colon(\C^N,0)\to (\C^N,0)$ sending $M$ into $M'$ which is transversally nonflat  satisfies 
  $\jac H \not \equiv 0$. 
\end{thm}

The proof of Theorem~\ref{t:hndhyp} will be obtained by following the lines of the proof of Theorem~\ref{t:hypinc}.
Note however that, contrarily to Theorem~\ref{t:hypinc}, the map $H$ in Theorem~\ref{t:hndhyp} need not be CR-transversal as the basic example $M=\{(z,w)\in \C^2: {\rm Im}\, w=|zw|^2\}$, $M'=\mathbb{H}_1$ and $H(z,w)=(z,zw)$, shows.

\subsection{Proofs of the results}

\begin{proof}[Proof of Theorem~{\rm \ref{t:hypinc}}]
Choose normal coordinates $Z=(z,w)$ for $M$ so that $M$ is given by \eqref{e:coord} and similarly for $M'$ (we just  add a ``prime" to the corresponding objects for the target hypersurface) and write $H=(F,G)=(F^1,\ldots, F^{N-1},G)$ as in \S~\ref{s:beginning}. 
Assume, by contradiction, that $H$ is totally degenerate. We will show that this
forces the transversal component $G$ to vanish.
Since $H$ sends $M$ into $M'$, we have the power series identity 
\begin{equation}
  G(z,Q(z,\chi,0)) = Q^\prime \left(
  F(z,Q(z,\chi,0)),\bar F (\chi,0),0 \right)
  \label{e:bsecseg}
\end{equation}
that we differentiate with respect to $\chi$ to
obtain
\begin{multline} 
  \left( G_w (z,Q(z,\chi,0)) - 
  Q^\prime_{z^\prime} 
  \left( F(z,Q(z,\chi,0)),\bar F(\chi,0),0 \right)\cdot
  F_w (z,Q(z,\chi,0)\right) \,Q_{\chi} (z,\chi,0)\\ = 
  Q_{\chi^\prime}^{\prime} \left(F(z,Q(z,\chi,0)), \bar F (\chi,0),0\right) 
  \cdot \dopt{\bar F}{\chi} (\chi,0).
   \label{e:diff1}
\end{multline}
Taking derivatives (of arbitrary order) with respect to $z$ and evaluating at $0$ we see
that the right hand side of \eqref{e:diff1} is always contained in the ($\C$-linear) span of the
$\bar F^j_{\chi} (\chi,0)$, $j=1,\ldots,N-1$; this space has by assumption 
dimension less than $N-1$ (strictly). If the coefficient of $Q_\chi (z,\chi,0)$
on the left hand side of \eqref{e:diff1} is not zero, 
an application of Lemma~\ref{l:prolongation} (with $b(z,\chi)=Q_\chi (z,\chi,0)$) yields that for every $\alpha \in \N^{N-1}$, 
$Q_{z^\alpha,\chi} (0,\chi,0)$ is contained in the (${\rm Frac}\, \C \dbl \chi\dbr$-linear) span of the 
$\bar F^j_{\chi} (\chi,0)$, $j=1,\ldots,N-1$, 
a contradiction, since $M\in {\mathcal C}$ means that
the dimension of the space spanned (over ${\rm Frac}\, \C \dbl \chi\dbr$) by the $Q_{z^\alpha,\chi} (0,\chi,0)$, $\alpha \in \N^{N-1}$, is $N-1$. 

Thus, we conclude that 
\begin{equation}
  G_w (z,Q(z,\chi,0)) - 
  Q^\prime_{z^\prime} 
  \left( F(z,Q(z,\chi,0)),\bar F(\chi,0),0 \right)\cdot
  F_w (z,Q(z,\chi,0)) = 0.
  \label{e:diff1concl}
\end{equation}
Setting $\chi =0 $ in this equation we have that 
$G_w (z,0) = 0$. We now construct by induction for each positive integer $k$
formal power series $c^k (z,w,\chi^\prime)$
satisfying
\begin{align}
  c^k (z,w,0) &= 0, \label{e:indhyp1}\\
  G_{w^k} (z,Q(z,\chi,0)) - c^k (z,Q(z,\chi,0) ,\bar f(\chi,0)) & = 0.
  \label{e:indhyp2}
\end{align}
In view of \eqref{e:diff1concl}, we have already constructed $c^1$. For every integer $\ell \geq 2$, if $c^1,\ldots,c^{\ell -1}$ have been constructed, we then   
differentiate the equation \eqref{e:indhyp2} 
for $k = \ell -1$ with respect to 
$\chi$ and we obtain 
\begin{multline}
  \left( G_{w^\ell} (z,Q(z,\chi,0)) - 
  c^{\ell-1}_{w} (z,Q(z,\chi,0),\bar F(\chi,0))
  \right) \, Q_\chi (z,\chi,0) \\= c^{\ell-1}_{\chi^\prime} (z,Q(z,\chi,0),
  \bar F(\chi,0)) \cdot \dopt{\bar F}{\chi} (\chi,0).
  \label{e:indstep1}
\end{multline}
The same argument as in the case $\ell=1$ above lets us conclude that
the coefficient of
$Q_\chi (z,\chi,0)$ on the left hand side of \eqref{e:indstep1}
vanishes. Thus, we can choose 
\[c^\ell(z,w,\chi^\prime) := c^{\ell-1}_{w} (z,w,\chi^\prime),\]
which satisfies \eqref{e:indhyp1} and \eqref{e:indhyp2} 
for $k = \ell$.  We may now conclude that 
$G_{w^\ell} (z,0) = 0$ for all $\ell$, so that the 
transversal component of $H$ is flat. The proof of the theorem is complete.
\end{proof}

For the proof of Corollary~\ref{c:newcor}, we need the following result contained in the proof of \cite[Proposition 7.1]{BER5}.

\begin{lem}\label{l:ber} Let $M,M'\subset \C^N$ be formal real hypersurfaces. Assume that
 $M$ is holomorphically nondegenerate and that $M'$ does not contain any formal curve. Then every formal holomorphic map sending $M$ into $M'$ that is transversally flat is necessarily constant.
\end{lem}

\begin{proof}[Proof of Corollary~{\rm \ref{c:newcor}}]
The corollary is obtained by combining Theorem~\ref{t:hypinc} and Lemma~\ref{l:ber} and noticing that if $M\in {\mathcal C}$ then $M$ is holomorphically nondegenerate.
\end{proof}

\begin{proof}[Proof of Proposition~{\rm \ref{c:add}}] The implication (i)$\Rightarrow $(ii) is proved in \cite[Corollary 4.2]{LM3} while the implication (i)$\Rightarrow$(ii) (with the weaker assumption that $M$ is of finite type instead of being in the class ${\mathcal C}$) is the content of \cite[Theorem 4.1]{ER1}.
\end{proof}

\begin{proof}[Proof of Theorem~{\rm \ref{t:hndhyp}}] We keep the notation of the proof of Theorem~\ref{t:hypinc} and present a completely analogous argument. We also write $\zeta=(\chi,\tau)\in \C^{N-1}\times \C$ and $H=(F,G)=(H^1,\ldots,H^N)$. Here again assume, by contradiction, that $\jac H  \equiv 0$. We will show that this
forces the transversal component $G$ to vanish.
Since $H$ sends $M$ into $M'$, we have the formal power series relation  
\begin{equation}
  G(z,Q(z,\zeta)) = Q^\prime \left(
  F(z,Q(z,\zeta)),\bar H (\zeta)\right).
  \label{e:bsecsegbis}
\end{equation}
We differentiate \eqref{e:bsecsegbis} with respect to $\zeta$ to
get
\begin{multline} 
  \left( G_w (z,Q(z,\zeta)) - 
  Q^\prime_{z^\prime} 
  \left( F(z,Q(z,\zeta)),\bar H(\zeta) \right)\cdot
  F_w (z,Q(z,\zeta)\right) \,Q_{\zeta} (z,\zeta)\\ = 
  Q_{\zeta^\prime}^{\prime} \left(F(z,Q(z,\zeta)), \bar H (\zeta)\right) 
  \cdot \dopt{\bar H}{\zeta} (\zeta).
   \label{e:diff1bis}
\end{multline}
Taking derivatives (of arbitrary order) with respect to $z$ and evaluating at $0$ we see
that the right hand side of \eqref{e:diff1bis} is always contained in the ($\C$-linear) span of the
$\bar H^j_{\zeta} (\zeta)$, $j=1,\ldots,N$; this space has by assumption 
dimension less than $N$ (strictly). If the coefficient of $Q_\zeta (z,\zeta)$
on the left hand side of \eqref{e:diff1bis} is not zero, 
another application of Lemma~\ref{l:prolongation} yields that for every $\alpha \in \N^{N-1}$, 
$Q_{z^\alpha,\zeta} (0,\zeta)$ is contained in the (${\rm Frac}\, \C \dbl \zeta\dbr$-linear) span of the 
$\bar H^j_{\zeta} (\zeta)$, $j=1,\ldots,N$, 
a contradiction, since $M$ being holomorphically nondegenerate means that
the dimension of the space spanned (over ${\rm Frac}\, \C \dbl \zeta\dbr$) by the $Q_{z^\alpha,\zeta} (0,\zeta)$, $\alpha \in \N^{N-1}$, is $N$ (see e.g.\ \cite{BERbook, S2}). 

We therefore conclude that 
\begin{equation}
  G_w (z,Q(z,\zeta)) - 
  Q^\prime_{z^\prime} 
  \left( F(z,Q(z,\zeta)),\bar H(\zeta) \right)\cdot
  F_w (z,Q(z,\zeta)) = 0.
  \label{e:diff1conclbis}
\end{equation}
Setting $\zeta =0 $ in \eqref{e:diff1conclbis} we get that 
$G_w (z,0) = 0$. As in the proof of Theorem~\ref{t:hypinc}, we define by induction for each positive integer $k$
formal power series $S^k (z,w,\zeta^\prime)$
satisfying
\begin{align}
  S^k (z,w,0) &= 0, \label{e:indhyp1bis}\\
  G_{w^k} (z,Q(z,\zeta)) - S^k (z,Q(z,\zeta) ,\bar H(\zeta)) & = 0.
  \label{e:indhyp2bis}
\end{align}
We have already constructed $S^1$ above. For every integer $\ell \geq 2$, if $S^1,\ldots,S^{\ell -1}$ have been defined, we then   
differentiate the equation \eqref{e:indhyp2bis} 
for $k = \ell -1$ with respect to 
$\zeta$ and obtain 
\begin{equation}
  \left( G_{w^\ell} (z,Q(z,\zeta)) - 
  S^{\ell-1}_{w} (z,Q(z,\zeta),\bar H(\zeta))
  \right) \, Q_\zeta (z,\zeta) = S^{\ell-1}_{\zeta^\prime} (z,Q(z,\zeta),
  \bar H(\zeta)) \cdot \dopt{\bar H}{\zeta} (\zeta).
  \label{e:indstep1bis}
\end{equation}
The same argument as in the case $\ell=1$ above lets us conclude that
the coefficient of
$Q_\zeta (z,\zeta)$ on the left hand side of \eqref{e:indstep1bis}
vanishes. Thus, we can choose 
\[S^\ell(z,w,\zeta^\prime) := S^{\ell-1}_{w} (z,w,\zeta^\prime),\]
which satisfies \eqref{e:indhyp1bis} and \eqref{e:indhyp2bis} 
for $k = \ell$.  All this implies that
$G_{w^\ell} (z,0) = 0$ for all $\ell$, so that the 
transversal component of $H$ is flat. The proof of the theorem is finished.
\end{proof}

%\begin{rem}
%In the case where $M,M'$ are real-analytic hypersurfaces through the origin and $H$ is a local holomorphic map sending %the germ $(M,0)$ into $(M',0)$, the proof of Theorem~\ref{t:hndhyp} is rather obvious and proceeds as follows: if $H$ %is transversally non-flat, there is a neighbourhood $\Sigma$ of $0$ in $M$ and a Zariski open subset $\Omega$ of %$\Sigma$ such for all $p\in \Omega$, $H$ is CR-transversal at $p$ and $M$ is finitely nondegenerate at $p$ (see e.g.\ %\cite{BERbook}). Then by \cite[Corollary 4.6]{LM3}, $H$ is a local biholomorphism near all such points $p$ and hence %$\jac H \not \equiv 0$.
%\end{rem}

\section{The infinite type case}\label{s:infinitetype}

In this section, we will derive some analogues of the transversality result due
to Baouendi and Rothschild \cite{BR4} for infinite type hypersurfaces.  We say
that a formal real hypersurface $M\subset \C^N$ given in normal coordinates
$(z,w)\in \C^N$ is of $m$-{\em infinite type}, where $m\in \Z_+$, if $M$ is
given by a (complexified) defining equation of the form \begin{equation} w=
  Q(z,\chi,\tau), \ \text{ where } \quad Q(z,\chi, \tau) = \tau + \tau^m
  \widetilde Q (z,\chi, \tau), \label{e:minftdef} \end{equation} with
  $\widetilde Q(z,0,\tau)=\widetilde Q(0,\chi,\tau)= 0$ and $$\widetilde Q (z,
  \chi, 0)\not \equiv 0.$$ (As customary, we use a similar notation for any
  other formal real hypersurface $M'$ by just adding a ``prime".) This above
  condition was introduced by Meylan \cite{Meylan95}, who also observed that
  such a condition together with the integer $m$ are invariantly attached to
  $M$ (i.e.\ independent of the chosen normal coordinates); the notion of 
  ``transversal order'' that we define next is also due to Meylan, who  
  established its basic properties.

\subsection{The transversal order of a map} Let $M,M'\subset \C^N$ be two
formal real hypersurfaces (given in normal coordinates) that we assume to be of
$m$-infinite and $m'$-infinite type respectively, $m,m'\geq 1$. Given a
transversally nonflat formal holomorphic map $H\colon (\C^N,0)\to  (\C^N,0)$
sending $M$ into $M'$ written $H = (F,G)$ as in \S \ref{s:beginning}, we define
the {\em transversal order} of $H$ and denote it by $\trord H$ to be the
integer $k$ such that $G (z,w) = w^k \widetilde{G} (z,w)$ with $\widetilde{G}
(z,0) \not \equiv 0$.  If $H$ is transversally flat, we just naturally set
$\trord H = \infty $. We note that in any case we always have $\trord H\geq 1$
and that the notion of transversal order is invariant under changes of coordinates preserving normal coordinates 
in the source and in the target space \cite[Lemma 2.1]{Meylan95}.

The first lemma we state
asserts that in the above situation we necessarily have $\widetilde{G}(0) \neq 0$. It is also contained in \cite{Meylan95} but for the reader's convenience we give a short proof.
\begin{lem}
  \label{l:trordnv} Let $M, M'$ be as above and  $H\colon (\C^N,0)\to (\C^N,0)$ be a formal holomorphic map sending $M$ into $M'$. If $H$ is 
  not transversally flat, that is, $\trord H = k < \infty$, then 
  $G_{w^k} (z,0)= G_{w^k} (0) \in \R\setminus\left\{ 0 \right\}$.
\end{lem}
\begin{proof}
  We start with the equation 
  \begin{equation}\label{e:basicFG}
    G\left( z,Q\left( z,\chi,\tau \right) \right) = Q^\prime 
    \left( F\left( z,Q\left( z,\chi,\tau \right) \right),\bar F \left( \chi,\tau \right),
    \bar{G}\left( \chi,\tau \right)\right),
  \end{equation}
  in which we set $\chi = 0$, $\tau = w$. This gives
  \begin{equation}\label{e:newstuff}
G(z,w) = \bar G (0,w) + \bar G(0,w)^{m'}
  \widetilde{Q}' (F(z,w),\bar F(0,w), \bar G(0,w)).
  \end{equation}
  The last equation clearly implies that $G(0,w)\neq 0$; so let $s$ be the
  smallest integer such that $G_{w^s}(0) \neq 0$. From the definition of $k$ and $s$, we immediately get
  that $s\geq k$. Furthermore we notice that  the order (in $w$) of left-hand side of \eqref{e:newstuff}
  is exactly $k$, whereas it is at least $s$ for the right-hand side of \eqref{e:newstuff} which shows that $k\geq s$
  and therefore $k=s$. Now comparing the coefficients of $w^k$ on both 
  sides of the equation we get $G_{w^k}(z,0) = \overline{G_{w^k}} (0)$, which completes the proof of the lemma.
\end{proof}
Our goal in this section is to bound the transversal order of a map in terms
 of $m$ and $m'$. The bound is provided in our next proposition.
 \begin{prop}
  \label{p:trordest} Let $M, M'$ be as above and  $H\colon (\C^N,0)\to (\C^N,0)$ be a formal holomorphic map sending $M$ into $M'$. If $H$ is 
  transversally nonflat, then
  \begin{equation}
    \left( m' - 1 \right)\, \trord H  \leq m - 1.
    \label{e:trordest}
  \end{equation}
\end{prop}

We will give examples illustrating this bound with specific maps
in Example~\ref{exa:blowups}.
Let us remark that if $m' = 1$, there are examples to show that no bound on 
the transversal order can be obtained in this case, see Example~\ref{exa:type1} below. Let us also add that if, in Proposition~\ref{p:trordest}, the map $H$ 
is assumed to be furthermore not totally degenerate, Meylan proved that the inequality given by \eqref{e:trordest} is in fact an equality (see \cite[Proposition 2.2]{Meylan95}).

\begin{proof}
  We again start with \eqref{e:basicFG}, set $\trord H =k$, and rewrite
  this equation in the following way:
  \begin{multline}\label{e:series}
  (\tau + \tau^m \tilde Q(z,\chi,\tau))^{k} \widetilde{G}\left( z, Q(z,\chi,\tau) \right)
  =\\
   \tau^k \overline{\widetilde{ G}} \left( \chi,\tau \right) + 
  \tau^{k m' } \overline{\widetilde{ G}} \left( \chi,\tau \right)^{m'}
  \widetilde{Q}' (F(z,Q(z,\chi,\tau)), \bar F(\chi,\tau) ,\bar G(\chi,\tau)).
  \end{multline}
  On the left hand side of this equation, we claim that there is 
  a term which has
  order  $k + m - 1$ in $\tau$ whose coefficient depends on $z$ (and $\chi$). Indeed,
  let us assume that $$\tilde Q (z,\chi,0) = \sum_{(a,b)\in (\Z_+)^2} \tilde Q_{a,b} (z,\chi,0),\quad 
  {\rm where}\quad  \tilde Q_{a,b} (\lambda z, \mu \chi,0) = \lambda^a \mu^b \tilde Q_{a,b} (z,\chi,0),\ (\lambda,\mu)\in \C^2,$$
  and let $(a_0,b_0)$ be minimal under the condition that $\tilde Q_{a,b} 
  (z,\chi,0) = 0$ for $a+ b < a_0 + b_0$ and $a< a_0$. We expand
  the series on the left hand side of \eqref{e:series} in terms $T_{a,b,c}(z,\chi,\tau)$, 
  which are homogeneous in $z$, $\chi$, and $\tau$ i.e.\ $T_{a,b,c} (\lambda z, \mu \chi, \delta \tau)
  = \lambda^a \mu^b \delta^c T_{a,b,c} (z,\chi,\tau)$ for $\lambda,\mu$ as above and $\delta \in \C$. We claim that
  $T_{a_0,b_0,k+m -1}(z,\chi,\tau) = k \widetilde{G} (0)\tau^{k+m-1} \tilde Q_{a_0,b_0} (z,\chi,0) \neq 0$.
  First note that modulo terms $T_{a,b,c}$ with $c > k+ m -1$, we can 
  rewrite the series as $(\tau + \tau^m \tilde Q (z,\chi,0))^k \widetilde{G}
  (z,\tau)$. We now pick the terms which are homogeneous of order $a_0$ in $z$ and
  $b_0$ in $\chi$ and see that our claim holds.
  
  Now, the lowest order in $\tau$ of a term of the series 
  on the right hand side of \eqref{e:series} which contains a coefficient depending on $z$
  is at least $km'$. Thus, $km' \leq k + m - 1$, 
  which is the inequality we wanted to prove.
\end{proof}

Given a formal real hypersurface $M\subset \C^N$, $M$ is of infinite type if and only if there
is a formal complex hypersurface $E$, which we refer to as the ``exceptional'' 
complex hypersurface, which is contained in $M$. (In normal coordinates, this
is necessarily given by $E = \left\{ w= 0 \right\}$.) Note that if $M$ is such a hypersurface, a necessary condition for a formal holomorphic self-map of $M$ to be CR-transversal is to be transversally non-flat i.e.\ $H(\C^N)\not \subset E$. As a direct application of Proposition~\ref{p:trordest} we obtain that  the converse of this statement does hold for $m$-infinite type hypersurfaces provided that $m\geq 2$. In addition, Example~\ref{exa:type1} below shows that the corresponding result does not hold for $1$-infinite type hypersurfaces. 

\begin{cor}\label{c:transtype2} Let $M\subset \C^N$
  be a formal real hypersurface of $m$-infinite type with 
  $m\geq 2$, and denote by $E$ the exceptional
  complex hypersurface contained in $M$. Then for every formal holomorphic map $H\colon (\C^N,0)\to (\C^N,0)$ sending $M$ into itself, either $H$ is CR-transversal or $H(\CN)\subset E$.
\end{cor}
Another application of Proposition~\ref{p:trordest} is given by the following generalization of Corollary~\ref{c:transtype2}.

\begin{cor}
  \label{c:nomaps}  Let $M,M'\subset \C^N$
  be formal real hypersurfaces of $m$-infinite and $m'$-infinite type respectively, and denote by $E'$ 
  the exceptional complex hypersurface contained in $M'$. 
  Assume that $1 < m'\leq m$ and that $m < 2m' - 1$. Then for every formal holomorphic map $H\colon (\C^N,0)\to (\C^N,0)$ sending $M$ into $M'$, either $H$ is CR-transversal or $H(\CN)\subset E'$.
\end{cor}
We will now discuss the examples we already pointed out.  Our first example provides explicit maps between 
certain infinite type hypersurfaces, where we can compute the infinite type and the 
transversal order of the maps explicitly.
\begin{exa}\label{exa:blowups}
  Let us note that (for different reasons) certain
  special cases of this example have 
  been considered by Kowalski \cite{Travis1}
  and Zaitsev \cite{Zsurvey}.

  Consider the Heisenberg hypersurface $\mathbb{H}_1\subset \C^2_{(z,w)}$ as given in Example~\ref{ex:example1}. We let 
  \[ H_{b,c} (z,w) = ( \sqrt{c} z w^b, w^c),\quad b,c\in \Z_+,  \]
  be a weighted analytic blowup of $\C^2$; for reasons
  which will become apparent below, we concentrate on blowups where 
  $2b>c$.
  The preimage of $\mathbb{H}_1$ under this weighted blowup 
  is given by 
  \[ \imag w^c = c |z|^{2} |w|^{2b}; \]
  setting $w = s + it$, this equation turns into
  \[ 
  \sum_{k=0}^{[(c-1)/2]} \binom{c}{2k+1} (-1)^k s^{c-2k-1} t^{2k+1}
  = c |z|^{2} \sum_{k=0}^{b} \binom{b}{k} s^{2b-2k} t^{2k},
  \]
  where $[\cdot]$ denotes the integer part.  We set $s^d u = t$, where $d=2b-c+1$; after 
  division by $s^{2b}$ the equation above now becomes
  \[
  \sum_{k=0}^{[(c-1)/2]} \binom{c}{2k+1} (-1)^k s^{2k (d-1)} {u}^{2k+1}
  = c |z|^{2} \sum_{k=0}^{b} \binom{b}{k} s^{2k (d-1)} {u}^{2k}.
  \]
  By the Implicit Function Theorem, this equation has a unique analytic solution 
  $u = \Theta_{b,c} (|z|^{2},s)$ which satisfies $\Theta_{b,c}(|z|^{2},0)
  =|z|^{2} \neq 0$ and $\Theta_{b,c}(0,s)=0$ (here we use that $2b > c$). Thus, the preimage of 
  $\mathbb{H}_1$ under the blowup $H_{b,c}$ contains a germ of 
  a real-analytic hypersurface
  $M_{b,c}$, of $d$-infinite type  given by the equation
  \[
  t = s^d \Theta_{b,c} (|z|^{2}, s).
  \]
  From this construction, we obtain examples of maps between 
  hypersurfaces of infinite type by noticing that 
  $H_{b,c}\circ H_{\tilde b,\tilde c}
  = H_{\tilde b + b \tilde  c, c \tilde c}$ with $2\tilde b>\tilde c$, $\tilde b,\tilde c\in \Z_+$; we claim that $H_{\tilde b,\tilde c}$ is a map of 
  transversal order $\tilde c$ from the $m$-infinite 
  type hypersurface $M_{\tilde b + b \tilde c,c\tilde c}$ to the
  $m'$-infinite type hypersurface 
  $M_{b,c}$ where 
  $m=2(\tilde b+b\tilde c)-c\tilde c+1$ and $m'=2b-c+1$; 
  for this map, the bound \eqref{e:trordest} becomes
  \[ 
  \trord H_{\tilde b, \tilde c} = \tilde c \leq \frac{m-1}{m'-1} 
  = \tilde c +\frac{2 \tilde b}{2b -c}.
  \]
  Let us establish our claim above i.e.\ that 
  $H_{\tilde b,\tilde c}$ actually sends 
  $M_{\tilde b + b \tilde c,c\tilde c}$ 
  into $M_{b,c}$. We rewrite the defining equation 
  of $M_{b,c}$ in complex form and replace the variables by 
  $(\zeta,\eta)$;
  then, $M_{b,c}\subset \C_{\zeta,\eta}^2$ is given by 
  \[ \eta = \bar \eta + \bar \eta^{m'} \Phi (|\zeta|^2, \bar \eta),\]
  where now $\Phi (|\zeta|^2,0) = 2 i |\zeta|^2$ and $\Phi (0,s)=0$. To prove our claim,
  we will show that there is exactly one germ of a real-analytic 
  hypersurface in $\C^2_{z,w}$
  of the form $ t = s^{m} \varphi (|z|^2, s)$, $w=s+it$, with $\varphi (|z|^2,0) = 
  |z|^2$ in the preimage of $M_{b,c}$ under the blowup
  $H_{\tilde b, \tilde c}$ with 
  $m = (m'-1) \tilde c  + 2\tilde b +1$; 
  since the preimage of $\mathbb{H}$ under 
  the composition $H_{b,c} \circ H_{\tilde b, \tilde c}$ is another such,
  the real-analytic hypersurface constructed in this way in the
  preimage of $M_{b,c}$ under $H_{\tilde b, \tilde c}$ agrees with 
  $M_{\tilde b + b \tilde c, c \tilde c}$. 

  As before, we substitute $\eta = w^{\tilde c}$, and 
  $\zeta = \sqrt{\tilde c} z w^{\tilde b}$.
  The equation of $H_{\tilde b,\tilde c}^{-1}(M_{b,c})$ 
  turns into
  \begin{equation}\label{e:sickofit}
 2i \sum_j \binom{\tilde c}{2j+1} 
  (-1)^{j} s^{\tilde c - (2j +1)} t^{2j+1}
  = (s- it)^{\tilde c m'} 
  \Phi ( \tilde c |z|^2 (s^2 + t^2)^{\tilde b}, (s - it)^{\tilde c}).
\end{equation}
  We now  proceed as above and substitute $t = \tilde t s^{m}$,
  where $m = (m'-1) \tilde c +1 + 2\tilde b$.
  Writing $\Phi$ as a Taylor series $ \Phi( |\zeta|^2 , \bar \eta )=\sum_{k\geq 1,l\geq 0} \Phi_{k,l} |\zeta|^{2k} \bar \eta^l$
  \[\begin{aligned} 
  \Phi( \tilde c |z|^2 (s^2 + t^2)^{\tilde b}, (s - it)^{\tilde c} ) &= 
  \sum_{k\geq 1,l\geq 0} \Phi_{k,l} \,
  (s^2+ \tilde t^2 s^{2m})^{\tilde b k} {\tilde c}^k | z|^{2k}(s-i\tilde ts^m)^{\tilde c l}\\
  &=
  s^{2\tilde b} \sum_{k\geq 1,l\geq 0} \Phi_{k,l} \,
  s^{2\tilde b (k-1)+\tilde c l} (1+ \tilde t^2 s^{2m-2})^{\tilde b k} {\tilde c}^k
  |z|^{2k}(1-i\tilde t s^{m-1})^{\tilde c l}\\
  &= s^{2 \tilde b}  2i \tilde c |z|^2 + O(s^{2 \tilde b + 1} )
\end{aligned}
  \]
  we see that the right hand side of \eqref{e:sickofit} 
  is divisible by $s^{m'\tilde c + 2\tilde b}$ (as is 
  the left hand side), and after division, 
  our equation becomes
  \[ 2i \tilde t
  =
  2i |z|^2 + O(s).\]
  Here, we have again tacitly used the fact that since $2b > c$ and $m'>1$, 
  $m > 1$, and 
  the Implicit Function Theorem guarantees a unique solution 
  $\tilde t = \varphi (|z|^2, s)$ with $\varphi (|z|^2,0)
  = |z|^2$ 
  of this equation as claimed.
\end{exa}
\begin{exa}
  \label{exa:type1} For each positive integer $k$, the map $T_k\colon \C^2\to \C^2$ given by $T_k(z,w) = (z, w^k)$ maps
   the $1$-infinite type hypersurface $M_k$ given by 
   \[ w = \bar w e^{\frac{i|z|^2}{k}} \]
   into the $1$-infinite type hypersurface $M$ given by
   \[ w = \bar w e^{i |z|^2}.\]
   Thus, the maps 
   \[
   H_k (z,w) = (\sqrt{k} z, w^k) 
   \]
   are examples of holomorphic self-maps sending the $1$-infinite type hypersurface $M$ into itself 
   having arbitrary transversal order.
\end{exa}

\subsection{On being an automorphism in the infinite type case}
We conclude this note by proving an analogue in the infinite type case 
of a criterion (established in \cite{LM3}) for a map to be an automorphism. We know from Corollary~\ref{c:transtype2} that if $M$ a formal real hypersurface that is of $m$-infinite type with $m\geq 2$, then every formal holomorphic self-map $H$ of $M$ is either
CR-transversal or transversally flat. We
now put an additional nondegeneracy assumption on $M$ that will force any formal CR-transversal self-map of $M$ to be  an automorphism (this will also be true in the 
$1$-infinite type case, where in view of Example
\ref{exa:type1}, the transversality result does not hold in general, see Theorem~\ref{p:inftypetvimpaut}). For this, we say that 
$M$ belongs to the class $\mathcal{C}_m$ if $M$ is of $m$-infinite type and is given in normal coordinates as follows
\[ w= Q(z,\chi,\tau), \quad Q(z,\chi,\tau) = \tau + \tau^m \widetilde{Q} (z,\chi,\tau), \]
with  $\widetilde{Q} (z,\chi,0) = \sum_{\alpha \in \N^{N-1}} \widetilde{Q}_\alpha (\chi) z^\alpha$, and with 
the formal map 
\begin{equation}
  \chi \mapsto \left( \widetilde{Q}_\alpha (\chi)  \right)_{|\alpha| \leq k}
  \label{e:isgenfin}
\end{equation}
of generic rank $N-1$ for large enough $k$. We leave it to the reader to check that this condition is independent of the choice of normal coordinates.

Our goal in this section is to establish the following theorem.

\begin{thm}
  \label{t:inftypeauto} Let $M\subset \C^N$ be a formal real hypersurface that belongs to the class ${\mathcal C}_m$ for some $m\geq 2$ and denote
   by $E$ the exceptional complex hypersurface contained in $M$. Then for every formal holomorphic map $H\colon (\C^N,0)\to (\C^N,0)$ sending $M$ into itself, one of the  following two conditions holds:
  \begin{compactenum}[i)]
  \item[{\rm (i)}] $H(\CN)\subset E$;
  \item[{\rm (ii)}] $H$ is an automorphism.
  \end{compactenum}
\end{thm}

Note here again that Example~\ref{exa:type1} shows that Theorem~\ref{t:inftypeauto} does not hold in the $1$-infinite type case.

We  start with the following lemma that can be found in \cite[Lemma 3.8]{Kow2} for which we provide 
a different and simpler proof.
\begin{lem}
  \label{l:basidinft} 
  Let $M\subset \C^N$ be a formal real hypersurface of $m$-infinite type with $m\geq 1$ and given in normal coordinates as above. Then for every formal CR-transversal holomorphic map $H\colon (\C^N,0)\to (\C^N,0)$ sending $M$ into itself, $H=(F,G)$, we have
  \begin{equation}
    \widetilde{Q} (z,\chi,0) = G_w (0)^{m-1} \widetilde{Q} \left( F(z,0),\bar F(\chi,0),0 \right).
    \label{e:basidinft}
  \end{equation}
\end{lem}
\begin{proof}
  For the proof, we start (as usual) with the identity 
  \begin{equation}\label{e:addit}
 G(z,Q(z,\chi,\tau)) = Q(F(z,Q(z,\chi,\tau)),\bar F(\chi,\tau), \bar G(\chi,\tau)),
 \end{equation}
and we write $G(z,w)=w\, \widetilde G(z,w)$ with $\widetilde G(z,0)=G_w(0)\in \R\setminus \{0\}$ in view of Lemma~\ref{l:trordnv}. We may rewrite \eqref{e:addit} as follows
\begin{multline}\label{e:chance}
 (\tau+\tau^m\widetilde Q(z,\chi,\tau))\widetilde G(z,Q(z,\chi,\tau))=\cr \bar G(\chi,\tau)+(\bar G(\chi,\tau))^m\,  \widetilde Q(F(z,Q(z,\chi,\tau)),\bar F(\chi,\tau),\bar G(\chi,\tau)).
\end{multline}
We look at the terms of order $m$ in $\tau$ in the above equation. On the left side, such a term is given by 
the quantity 
\begin{equation}\label{e:quantity1}
 \widetilde G_{w^{m-1}}(z,0)+\widetilde Q(z,\chi,0) G_w(0),
\end{equation}
while on the right hand side, it is given by
\begin{equation}\label{e:quantity2}
\bar G_{\tau^m}(\chi,0)+(\bar G_\tau(0))^m \widetilde Q(F(z,0),\bar F(\chi,0),0).
\end{equation}
Since \eqref{e:quantity1}=\eqref{e:quantity2}, we first get, after setting consecutively $z=0$ and $\chi=0$ in the obtained identity that $\widetilde G_{w^{m-1}}(z,0)$ and $\bar G_{\tau^m}(\chi,0)$ are equal and constant. This leads to the desired identity. The proof of the lemma is complete.
\end{proof}
To prove Theorem~\ref{t:inftypeauto}, we need the following general fact about power series.

\begin{lem}
  \label{l:easystuff} Let $A\colon (\Cn_z\times\Cn_\chi,(0,0)) \to \C$ be a formal power series such that for $k$
  large enough the map
  \[ \chi\mapsto \left( A_{z^\alpha} \left( 0,\chi \right) \right)_{|\alpha|\leq k}\]
  is of generic rank $n$. If for some formal map $B\colon (\C^n,0)\to (\C^n,0) $ there is a
  nonzero constant $r$ such that
  \[ A(z,\chi) = r A (B(z),\bar B(\chi)), \]
  then necessarily $B$ is a formal automorphism of $\Cn$.
\end{lem}
\begin{proof}
  From the equation $A(z,\chi) = r A (B(z), \bar B(\chi))$ we 
  obtain for each $\alpha\in\Nn$ a polynomial map $\Phi_\alpha$ such that  
  $A_{z^\alpha} (0,\chi) = \Phi_{\alpha} \left( 
  \left(A_{z^{\beta}} (0,\bar B(\chi))\right)_{|\beta|\leq |\alpha|}\right).$
  We choose $k$ large enough such that 
  \[ \chi\mapsto\left( A_{z^\alpha} (0,\chi) \right)_{|\alpha|\leq k} \]
  is of generic rank $n$ and also choose $n$ multiindeces $\alpha^{(1)},\dots,\alpha^{(n)}$ such
  that the order $\nu$ of the determinant of the matrix
  \[ \begin{pmatrix}
    \frac{ \partial A_{z^{\alpha^{(1)}}}}{\partial \chi_1} (0,\chi) & \dots& \frac{ \partial A_{z^{\alpha^{(1)}}}}{\partial \chi_n} (0,\chi) \\
    \vdots & &\vdots \\
    \frac{ \partial A_{z^{\alpha^{(n)}}}}{\partial \chi_1} (0,\chi) & \dots& \frac{ \partial A_{z^{\alpha^{(n)}}}}{\partial \chi_n} (0,\chi) 
  \end{pmatrix}\]
  is minimal among any such choices of $n$ multiindeces. Differentiating the $n$ equations
  \[A_{z^{\alpha^{(j)}}} (0,\chi) = 
  \Phi_{\alpha^{(j)}}\left( (A_{z^\beta} (0,\bar B(\chi))_{|\beta|\leq|\alpha^{(j)}| }
  \right), \quad j=1,\dots,n,\]
  with respect to $\chi$ and using the Cauchy-Binet formula yields the inequality $\nu \geq \nu + \ord \det  \bar B_{\chi} (\chi)$. Consequently, $\ord \det \bar B_\chi (\chi)=0$ i.e.\ $\det B_z (0) \neq 0$.
\end{proof}
Combining the last two lemmas, we obtain the following result whose analogous statement 
 is contained in \cite{LM3} for hypersurfaces belonging to the class ${\mathcal C}$. 
\begin{thm}
  \label{p:inftypetvimpaut} Let $M\subset \C^N$ be a formal real hypersurface that belongs to the class ${\mathcal C}_m$ for some $m\geq 1$. Then every formal CR-transversal holomorphic map $H\colon (\C^N,0)\to (\C^N,0)$ sending $M$ into itself is an automorphism.
\end{thm}
\begin{proof}
  [Proof of Theorem~{\rm \ref{t:inftypeauto}}] The theorem is just obtained as the conjunction of Theorem~\ref{p:inftypetvimpaut} and Corollary~\ref{c:transtype2}.
\end{proof}

\bibliographystyle{plain}
\bibliography{bibliography_oct05}

\end{document}